\documentclass{article}

\usepackage{amsmath,amssymb,amsfonts}
\usepackage{mathrsfs}
\usepackage{graphicx}
\usepackage{color}
\usepackage{mathptmx}
\usepackage{setspace}
\usepackage{geometry}
\usepackage{indentfirst}
\usepackage{indentfirst}
\usepackage{url}
\usepackage{indentfirst}

\begin{document}

\begin{center}
\LARGE\bfseries Completely independent Steiner trees and corresponding tree connectivity

\vspace{1cm}

\normalsize
Jun Yuan\textsuperscript{a,*} \footnotemark, \footnotetext{Corresponding author. E-mail: junyuan1@tyust.edu.cn

{\textbf{Publisher's Note:} This is an Accepted Manuscript of an article published by Taylor \& Francis Group in the \emph{International Journal of Parallel, Emergent \& Distributed Systems} on \textbf{20/12/2025}, available online: \url{http://www.tandfonline.com/ } (DOI: 10.1080/17445760.2025.2609134).}
 }

Shan Liu\textsuperscript{a},
Shangwei Lin\textsuperscript{b},
Aixia Liu\textsuperscript{a}

\vspace{0.5cm}

\small
\textit{a. School of Applied Sciences, Taiyuan University of Science and Technology, Taiyuan, Shanxi 030024, People's Republic of China}

\textit{b. School of Mathematics and Statistics, Shanxi University, Taiyuan, Shanxi 030006, People's Republic of China}
\end{center}

\vspace{0.5cm}

\begin{center}
\rule{0.8\textwidth}{0.4pt}
\end{center}

\begin{center}
\textbf{Abstract}
\end{center}

\indent The $S$-Steiner tree packing problem provides mathematical foundations for optimizing multi-path information transmission, particularly in designing fault-tolerant parallelized routing architectures for massive-scale network infrastructures. In this article, we propose the definitions of completely independent $S$-Steiner trees (CISSTs for short) and generalized $k^*$-connectivity, which generalize the definitions of internally disjoint $S$-Steiner trees and generalized $k$-connectivity. Given a connected graph $G = (V,E)$ and a vertex subset $S\subseteq V, |S|\geq 2,$ an $S$-Steiner tree of $G$ is a subtree in $G$ that spans all nodes in $S.$ The $S$-Steiner trees $T_1,T_2,\cdots, T_k$ of $G$ are completely independent pairwise if for any $1\leq p<q\leq k,$ $E(T_p)\cap E(T_q)=\emptyset$, $V(T_p)\cap V(T_q)=S,$ and for any two vertices $x_{1},x_{2}$ in $S$, the paths connecting $x_{1}$ and $x_{2}$ in $T_p,T_q$ are pairwise internally disjoint. The packing number of CISSTs, denoted by $\kappa^*_G(S),$ is the maximum number of CISSTs in $G.$ The generalized $k^*$-connectivity $\kappa_k^*(G)$ is the minimum $\kappa_G^*(S)$ for $S$ ranges over all $k$-subsets of $V(G).$ We provide a detailed characterization of CISSTs. Also, we investigate the CISSTs of complete graphs and complete bipartite graphs. Furthermore, we determine the generalized $k^*$-connectivity for complete graphs and give a tight lower bound of the generalized $k^*$-connectivity for complete bipartite graphs.

\vspace{0.5cm}


\begin{center}

\end{center}

\section{Introduction}

Let $G$ be a finite undirected simple connected graph with vertex set $V(G)$ and edge set $E(G).$ For a vertex subset $S\subseteq V, |S|\geq 2,$ a subtree $T$ of $G$ is an $S$-Steiner tree if $T$ spans $S,$ i.e., all the leaves in $T$ belong to $S.$ Let $T_i(1\leq i\leq k)$ be $k(\geq 2)$ $S$-Steiner trees of a graph $G.$ Then $T_i(1\leq i\leq k)$ are called internally disjoint $S$-Steiner trees if for any $1\leq p<q\leq k,$ $E(T_p)\cap E(T_q)=\emptyset, V(T_p)\cap V(T_q)=S.$ We can refer to \cite{DuH} for applications of Steiner trees in computer communication networks.

Let $G$ be a connected graph. The connectivity $\kappa(G)$ of a graph $G$ is defined as the cardinality of the minimum vertex cut $F$ of $G$ provided $F$ exists; otherwise, $\kappa(G)$ is defined as $|V(G)|-1.$ Let $x_{1}$ and $x_{2}$ be a pair of distinct vertices in $G.$ Two paths from $x_{1}$ to $x_{2}$ are internally disjoint if they have no common edge and no common vertex except the terminal vertices $x_{1}$ and $x_{2}.$ The local connectivity of $x_{1},x_{2}$ in $G,$ denoted by $\kappa_G(x_{1},x_{2}),$ is defined as the maximum number of internally disjoint $(x_{1},x_{2})$-paths in $G.$ The Whitney Theorem \cite{Whitney2} shows that the connectivity $\kappa(G)=\min\{\kappa_G(x_{1},x_{2})\mid x_{1},x_{2}\in V(G), x_{1}\neq x_{2}\}.$ Similarly, for a vertex subset $S\subseteq V, |S|\geq 2,$ the local connectivity $\kappa_G(S)$ of $S$ in $G$ can be defined as the maximum number of internally disjoint $S$-Steiner trees in $G.$ The generalized $k$-connectivity $\kappa_k(G)$ \cite{Hager1985}, also known as $k$-set tree connectivity\cite{Hager1985}, is the minimum $\kappa_G(S)$ for $S$ ranges over all $k$-subsets of $V(G).$

The concept of generalized $k$-connectivity has attracted significant research interest. Extensive studies have been devoted to this property for various network structures, including hypercubes~\cite{Lin}, exchanged hypercubes~\cite{Zhao'}, dual cubes~\cite{Zhao2}, $(n,k)$-star networks~\cite{Li2020}, and hierarchical cubic networks~\cite{Zhao21}. Most of these efforts, however, are confined to the case where $3 \leq k \leq 4$. For arbitrary $k$, the generalized $k$-connectivity has been precisely determined for only a few graph classes, such as complete graphs~\cite{Chartrand}, complete bipartite graphs~\cite{Li2010'}, and complete equipartition $3$-partite graphs~\cite{Li2014}. Importantly, Li et al.~\cite{Li2012} established the computational complexity of the problem: deciding whether $\kappa_G(S) \geq l$ for a vertex subset $S$ (with $|S| \geq 2$) and an integer $l$ ($2 \leq l \leq |V(G)| - 2$) is NP-complete. They also provided tight bounds on $\kappa_3(G)$ for general graphs, as detailed in~\cite{Li2010}.

Let $T_1,T_2,\cdots,T_k$ be $k$ spanning trees of a graph $G,$ where $k\geq 2.$ Let $w$ be a vertex of $G.$ If for every vertex $v(\neq w)\in V(G)$, the $(w,v)$-paths in $T_i,1\leq i\leq k$ are pairwise internally disjoint, then we say that $T_1,T_2,\cdots,T_k$ are $k$ independent spanning trees rooted at $w.$ Numerous studies have presented constructions of independent spanning trees for a specified root vertex. However, reconstruction is unnecessary if a set of spanning trees remains independent under any choice of root vertices. Inspired by this perspective, Hasunuma \cite{Hasunuma2} proposed the definition of completely independent spanning trees. The $k$ spanning trees $T_i(1\leq i\leq k)$ are called completely independent spanning trees (CISTs for short) if they are edge-disjoint, and for any pair of vertices $x_{1}$ and $x_{2}$ in $G$, the paths connecting $x_{1}$ to $x_{2}$ in different trees are internally disjoint.

\begin{figure}[htbp]
  \centering
  \includegraphics[width=0.8\textwidth]{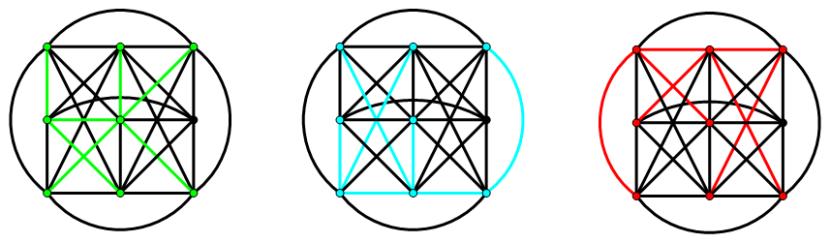}
  \caption{the $3$ CISSTs in graph $G$}
  \label{fig:figure1}
\end{figure}

Inspired by the definition of CISTs, the $k (\geq 2)$ $S$-Steiner trees $T_1,T_2,\cdots,T_k$ are called completely independent $S$-Steiner trees (CISSTs for short) if for any $1\leq p<q\leq k,$ $E(T_p)\cap E(T_q)=\emptyset$ and $V(T_p)\cap V(T_q)=S,$ for any $x_{1},x_{2} \in S$, their connecting paths in distinct $S$-Steiner trees $T_p,T_q$ are internally disjoint. For example, in a graph $G$ with $9$ vertices, there exist $3$ CISSTs with $|S|=8,$ see Figure 1. The packing number of CISSTs is the maximum cardinality of CISSTs in $G,$ denoted by $\kappa^*_G(S).$ The generalized $k^*$-connectivity $\kappa_k^*(G)$ is the minimum $\kappa_G^*(S)$ as $S$ ranges over all $k$-subsets of $V(G).$ By the definitions of internally disjoint $S$-Steiner trees and CISSTs, the CISSTs in $G$ are also internally disjoint $S$-Steiner trees; the converse is not necessarily true. Thus the generalized $k^*$-connectivity $\kappa_k^*(G)$ is at most the generalized $k$-connectivity $\kappa_k(G).$

In network reliability design, the CISTs and internally disjoint Steiner trees are essential theoretical tools for creating highly robust networks. The CISSTs become the CISTs when the set $S$ equals the entire vertex set $V(G)$ of graph $G$ (i.e., $S=V(G)$). The CISTs demand both edge-disjointness and internal vertex-disjointness among all spanning trees, providing dual fault tolerance against link and node failures. In contrast, CISSTs extend this independence to Steiner trees. Compared to internally disjoint Steiner trees, CISSTs impose much stricter independence requirements.

To illustrate the practical value of CISSTs, we introduce two practical scenario examples with the same terminal set
\[
S = \{u \text{ (Command Center)}, x \text{ (Hospital)}, v \text{ (Fire Station)}\}.
\]
In the emergency communication network $G = (V, E)$, the available relay nodes are $a$ and $b$. Figure 2 shows a traditional design based on two internally disjoint Steiner trees $T_1$ and $T_2$, where
\[
V(T_1) = \{u, a, x, b, v\},\quad E(T_1) = \{(u, a), (a, x), (x, b), (b, v)\},
\]
\[
V(T_2) = \{u, x, v\},\quad E(T_2) = \{(u, x), (x, v)\}.
\]
Although $T_1$ and $T_2$ are internally disjoint Steiner trees, the connecting path from $u$ to $v$ must pass through department $x$. This means that when department $x$ fails, the remaining network in Figure 2 no longer contains an $(S \setminus \{x\})$-Steiner tree. Consequently, the communication between $u$ and $v$ is disrupted and cannot be maintained.

\begin{figure}[htbp]
    \centering
    \begin{minipage}{0.32\textwidth}
        \centering
        \includegraphics[width=0.9\linewidth]{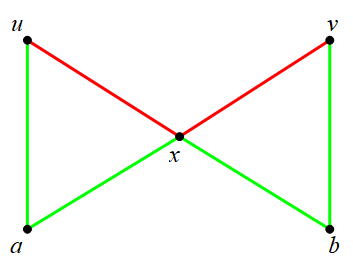}
        \caption{Traditional design with two internally disjoint Steiner trees}
        \label{fig:traditional}
    \end{minipage}
    \hspace{2cm}
    \begin{minipage}{0.32\textwidth}
        \centering
        \includegraphics[width=0.9\linewidth]{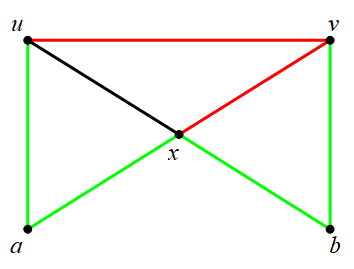}
        \caption{Enhanced design with two CISSTs}
        \label{fig:enhanced}
    \end{minipage}
\end{figure}

Now consider an enhanced network \( G' = (V', E') \), the available relay nodes be \( a, b \). As shown in Figure 3, we can construct two CISSTs \( T_1' \) and \( T_2' \), where,
\[
V(T_1') = \{u, a, x, b, v\},\quad E(T_1') = \{(u, a), (a, x), (x, b), (b, v)\},
\]
\[
V(T_2') = \{u, v ,x \},\quad E(T_2') = \{(u, v), (v, x)\}.
\]

The key advantage is that even when department $x$ fails, the system still has an $(S \setminus \{x\})$-Steiner tree, ensuring that departments $u$ and $v$ continue to operate effectively. This comparison demonstrates that CISSTs are practically essential for scenarios demanding guaranteed service continuity.

In fact, we derive a characterization of CISSTs: for any vertex $u\in S$, it can serve as an internal vertex in at most one $S$-Steiner tree in $G$. This property ensures that if there exist $k$ CISSTs in $G$ and a vertex $u\in S$ fails, then there must be $(k-1)$ completely independent $(S\setminus\{u\})$-Steiner trees in $G.$ However, when there exist $k$ internally disjoint $S$-Steiner trees in $G$ and a vertex $u\in S$ fails, it is not necessarily true that there are $(k-1)$ internally disjoint $(S\setminus\{u\})$-Steiner trees in $G.$ In this sense, compared with internally disjoint $S$-Steiner trees, CISSTs can naturally avoid shared bottlenecks and localize the effect of a single-point failure, thereby providing a provably stronger reliability guarantee for the design of ultimate fault-tolerant networks in fields such as finance, energy, and military systems.

From the perspective of computational complexity, determining the exact value of the generalized $k^*$-connectivity $\kappa_k^*(G)$ is generally a very difficult problem. This difficulty is closely related to the computational hardness of its classical counterpart---namely, determining the number of internally disjoint Steiner trees. Existing research has shown that: for any fixed integer $k \geq 4$, deciding whether a given $k$-subset $S \subseteq V(G)$ satisfies $\kappa_G(S) \geq \ell$ is NP-complete \cite{Li2012}; similarly, for any fixed $k \geq 3$, deciding whether $\kappa_k(G) \geq \ell$ is also NP-complete \cite{Chen L}.

Since a set of CISSTs is also a set of internally disjoint Steiner trees and requires that paths connecting the same terminals be internally vertex-disjoint, the problem of deciding whether $\kappa_G^*(S) \geq \ell$ imposes stronger conditions. Consequently, it is computationally at least as hard as the classical problems mentioned above, and thus is also NP-complete for $k \geq 4$.

This complexity motivates the search for exact polynomial-time algorithms on specific graph classes. The primary focus of this paper, however, is to conduct a thorough analysis of the generalized \( k^* \)-connectivity for two fundamental graphs: complete graphs and complete bipartite graphs. We determine the exact value of \( \kappa_k^*(K_n) \) and establish tight lower bounds for \( \kappa_k^*(K_{m,n}) \). Crucially, our proofs are constructive. Although a detailed algorithmic analysis beyond the scope of this work, the constructive nature of our proofs directly implies that the corresponding structures can be built in polynomial time for these graph classes. This work provides a basis for future explicit algorithmic formulations and analysis.

The remainder of this paper is organized as follows: in Section 2, we establish a theoretical characterization for CISSTs and investigate the properties of CISSTs and $\kappa_k^*(G)$ of $G.$ Additionally, we determine the generalized $k^*$-connectivity $\kappa_k^*(G)$ of complete graphs. In Section 3, we study the CISSTs in complete bipartite graphs and establish a tight lower bound for their $\kappa_k^*(G)$. Finally, our conclusions are presented in Section 4.

\section{A characterization of completely independent $S$-Steiner trees}

In this section, we systematically analyze the properties of the generalized connectivity parameter $\kappa_k^*(G)$ for general graphs $G.$ Subsequently, we derive an explicit formula for $\kappa_k^*(G)$ in the special case of complete graphs.

{\bf Theorem 2.1}. Let $T_i(1\leq i\leq k)$ be $k$ $S$-Steiner trees in $G.$ Then $T_i(1\leq i\leq k)$ are CISSTs if and only if they are edge-disjoint, and there is at most one $T_i$ for which $d_{T_i}(w)>1$ for any $w\in V(G).$

{\bf Proof.} Let $T_i(1 \leq i \leq k)$ be $k$ edge-disjoint $S$-Steiner trees, and suppose that for every vertex $w \in V(G)$, there is at most one $T_i$ such that $d_{T_i}(w) > 1$. Now, we show that for any $p \neq q$, $V(T_p) \cap V(T_q) = S$. Clearly $S \subseteq V(T_p) \cap V(T_q)$. Suppose, for a contradiction, that there exists a vertex $x' \in (V(T_p) \cap V(T_q)) \setminus S$. Since in an $S$-Steiner tree every leaf belongs to $S$, the vertex $x'$ is not a leaf in either $T_p$ or $T_q$. Hence $d_{T_p}(x') > 1$ and $d_{T_q}(x') > 1$, which contradicts the hypothesis. Thus $V(T_p) \cap V(T_q) = S$. Next, suppose there are two $S$-Steiner trees $T_p,T_q$ such that they are not completely independent, that is, in $T_p$ and $T_q$, a pair of vertices $x_{1},x_{2}\in S$ has paths that are not internally disjoint, say $(x_{1},x_{2})$-paths $P_{1}$ and $P_{2}$. Note that $T_p,T_q$ are edge-disjoint. Thus, there exists a common internal vertex $w$ in paths $P_{1}$ and $P_{2}$. It follows that $d_{T_p}(w)\geq 2$ and $d_{T_q}(w)\geq 2,$ contradicting the assumption that there exists at most one $T_i$ for which $d_{T_i}(w)>1.$

Let $T_i$ $(1 \leq i \leq k)$ be CISSTs. By definition they are edge-disjoint. We shall prove that for every vertex $w$, there is at most one $T_i$ with $d_{T_i}(w) > 1$. Suppose, for contradiction, that there exists a vertex $x$ and two distinct CISSTs, say $T_1$ and $T_2$, such that $d_{T_1}(x) \geq 2$ and $d_{T_2}(x) \geq 2$. Since $T_1$ and $T_2$ are CISSTs, we have $V(T_1) \cap V(T_2) = S$, whence $x \in S$. Since $d_{T_1}(x) \geq 2,$ there exists a $(a,b)$-path $P$ in $T_1$ such that $P$ passes through the vertex $x,$ and $a,b \in S\setminus\{x\}.$ Since $T_1$ and $T_2$ are two CISSTs, the $(a,b)$-path in $T_2$ must not pass through the vertex $x.$ And by $d_{T_2}(x) \geq 2,$ there exists a vertex $c\in S$ such that both the $(a,c)$-path and the $(b,c)$-path in $T_2$ pass through the vertex $x.$ Consider the $(a,c)$-path and $(b,c)$-path in $T_1$. Since $T_1$ and $T_2$ are two CISSTs, neither the $(a,c)$-path nor the $(b,c)$-path in $T_1$ pass through the vertex $x.$ It follow that the $(a,b)$-path in $T_1$ does not pass through the vertex $x,$ which contradicts the $(a,b)$-path $P$ in $T_1$ passes through $x.$
$\Box$

A vertex $v$ in a tree $T$ is said to be an internal vertex if $d_T(v) \geq 2$. Theorem 2.1 says that each vertex in graph $G$ is an internal vertex of at most one of the CISSTs. Next, we will discuss how the packing number of CISSTs changes as $S$ increases.

{\bf Theorem 2.2}. Let $S$ and $S'$ be two vertex subsets in $G$ with $S\subset S'.$ Then there are $k$ CISSTs in $G$ if there are $k$ completely independent $S'$-Steiner trees in $G,$ that is, $\kappa_G^*(S)\geq \kappa_G^*(S').$

{\bf Proof.} Assume $T_i'(1\leq i\leq k)$ are $k$ completely independent $S'$-Steiner trees in $G.$ For any $1\leq i\leq k,$ let $U_{i,1}=\{u\mid u\in V(T_i')\setminus S , d_{T_i'}(u)=1\},$ and let $T_{i,1}=T_i'-U_{i,1}.$ Let $U_{i,2}=\{u\mid u\in V(T_{i,1})\setminus S , d_{T_{i,1}}(u)=1\},$ and let $T_{i,2}=T_{i,1}-U_{i,2}.$ For convenience, denote $T_i'$ by $T_{i,0}.$ Repeat this process until for some $l\geq 1,$ $U_{i,l}$ is an empty set. Then we can obtain an $S$-Steiner tree $T_{i,l-1},$ and let $T_i=T_{i,l-1}.$

Next, we shall show the $S$-Steiner trees $T_i(1\leq i\leq k)$ are $k$ CISSTs in $G$. Since $T_i'(1\leq i\leq k)$ are completely independent, they are edge-disjoint, and hence $T_i(1\leq i\leq k)$ are also edge-disjoint. Since $T_i(1\leq i\leq k)$ are all $S$-Steiner trees, we have $S\subseteq V(T_p)\cap V(T_q)$ for any $1\leq p,q\leq k,p\neq q.$ Suppose there are two $S$-Steiner trees $T_p$ and $T_q$ such that $(V(T_p)\cap V(T_q))\setminus S\neq \emptyset,$ that is, there is a vertex $x_{1}\in (V(T_p)\cap V(T_q))\setminus S.$ From the construction of $T_i,$ we can deduce that $d_{T_{p}}(x_{1})>1$ and $d_{T_{q}}(x_{1})>1.$ It follows that $d_{T_{p}'}(x_{1})>1$ and $d_{T_{q}'}(x_{1})>1.$ On the other hand, since the $S'$-Steiner trees ${T_p}'$ and ${T_q}'$ are completely independent, by Theorem 2.1, we can deduce that at most one of $d_{T_{p}'}(x_{1})>1$ and $d_{T_{q}'}(x_{1})>1$ is true, a contradiction. Therefore, $S=V(T_p)\cap V(T_q)$ for any $1\leq p,q\leq k,p\neq q.$ Since $T_i'(1\leq i\leq k)$ are completely independent, for any two vertices $x_{2},x_{3} \in S,$ the $(x_{2},x_{3})$-paths in each $T_i'(1\leq i\leq k)$ are internally disjoint pairwise. And by the construction of $T_i(1\leq i\leq k)$, for any two vertices $x_{2},x_{3} \in S$, the $(x_{2},x_{3})$-paths in each $T_i, 1\leq i\leq k,$ are also internally disjoint pairwise. Therefore, $T_i(1\leq i\leq k)$ are $k$ CISSTs. According to the definition of $\kappa_G^*(S),$ we have $\kappa_G^*(S)\geq \kappa_G^*(S').$ The proof is complete.$\Box$

{\bf Corollary 2.3}. Let $k$ and $l$ be two integers with $2 \leq k\leq l.$ Then $\kappa_k^*(G)\geq \kappa_l^*(G).$

{\bf Proof.} Since $\kappa_k^*(G) = \min \big\{ \kappa_G^*(S) \mid S \subseteq V(G), |S| = k \big\},$ there exists a subset $S_{0}$ with $|S_{0}|=k$ such that $\kappa_G^*(S_{0}) = \kappa_k^*(G).$ By $\kappa_l^*(G) = \min \big\{ \kappa_G^*(S) \mid S \subseteq V(G), |S| = l\big\},$ for any subset $S_{1} \supseteq S_{0}$ with $|S_{1}|=l,$ it holds that $\kappa_G^*(S_{1}) \geq \kappa_l^*(G).$ By Theorem 2.2, we have $\kappa_G^*(S_{0}) \geq \kappa_G^*(S_{1}).$ Therefore,$\kappa_k^*(G) = \kappa_G^*(S_0) \geq \kappa_G^*(S_1) \geq \kappa_l^*(G).$ The proof is complete.$\Box$

{\bf Theorem 2.4}. Let $S$ be a vertex subset of order $k\geq 2$ in $V(G).$ Then $\kappa_k^*(G)\leq \kappa_k^*(G[S])+(|V(G)|-k).$

{\bf Proof.} Let $\mathcal{T}$ be a maximum set of CISSTs in $G.$ Let ${\mathcal{T}}_1$ be the set of trees in $\mathcal{T}$ whose vertex sets are all $S$, and let ${\mathcal{T}}_2= \mathcal{T} \setminus$ ${\mathcal{T}}_1.$ According to the definition of $S$-Steiner trees, all vertices in $V(G)\setminus S$ are internal vertices. Combined with Theorem 2.1, it follows that $|{\mathcal{T}}_2|\leq |V(G)|-k.$ Clearly, $|{\mathcal{T}}_1|\leq \kappa_k^*(G[S])$, and $\kappa_k^*(G)=|{\mathcal{T}}|=|{\mathcal{T}}_1|+|{\mathcal{T}}_2|.$ Thus $\kappa_k^*(G)\leq \kappa_k^*(G[S])+(|V(G)|-k).$ The proof is complete.$\Box$

{\bf Lemma 2.5} \cite{Pai}. There are $\lfloor\frac{n}{2}\rfloor$ CISTs in a complete graph $K_n$ for all $n\geq 4.$

{\bf Theorem 2.6}. Let $S$ be a vertex subset of order $s\geq 2$ in complete graph $K_n.$ Then there are $n-\lceil\frac{s}{2}\rceil$ CISSTs in $K_n.$

{\bf Proof.} By Lemma 2.5, there exist $\lfloor\frac{s}{2}\rfloor$ CISTs $T_1',T_2',\cdots, T_{\lfloor\frac{s}{2}\rfloor}'$ in the induced subgraph $K_n[S]$ (for convenience, we say there exists a CIST in the induced subgraph $K_n[S]$ when $s\leq 3.$). By the definitions of CISTs and CISSTs, $T_1',T_2',\cdots, T_{\lfloor\frac{s}{2}\rfloor}'$ are also the CISSTs in $K_n.$ By Theorem 2.1, for any $w\in S,$ there is at most one $T_i'$ such that $d_{T_i'}(w)>1.$

Next, we construct $n-s$ new $S$-Steiner trees $T_i.$ Denote $V(K_n)\setminus S=\{v_1,v_2,\cdots, v_{n-s}\}.$ Let $V(T_i)=S\cup\{v_i\},$ and let $E(T_i)=\{v_iu|u\in S\}.$ Clearly, for any $v\in V(K_n),$ there exists at most one $T_i$ such that $d_{T_i}(v)>1.$ From the construction of these $S$-Steiner trees, they are edge-disjoint. By Theorem 2.1, $T_1',T_2',\cdots, T_{\lfloor\frac{s}{2}\rfloor}'$ and $T_1,T_2,\cdots, T_{n-s}$ are $n-\lceil\frac{s}{2}\rceil$ CISSTs in $K_n.$ The proof is complete.$\Box$

Combining Theorem 2.4 and Theorem 2.6, we can deduce the following corollary.

{\bf Corollary 2.7}. Let $n$ and $s$ be two integers with $2 \leq s\leq n$ and $n\geq 4.$ Then the generalized $s^*$-connectivity $\kappa_s^*(K_n)=n-\lceil\frac{s}{2}\rceil.$

\section{Completely independent $S$-Steiner trees in complete bipartite graphs}

Let $(X,Y)$ be a bipartition of bipartite graph $K_{m_{1},m_{2}}$(with $2 \leq m_{1} \leq m_{2}$), where $X=\{x_1,x_2,\cdots, x_{m_{1}}\}, Y=\{y_1,y_2,\cdots, y_{m_{2}}\}.$ By the symmetry of $K_{m_{1},m_{2}},$ we assume $S=\{x_1,x_2,\cdots,x_i,y_1,y_2,\cdots,y_{s-i}\}.$ Clearly, $\max\{0,s-m_{2}\}\leq i\leq \min\{m_{1},s\}.$

{\bf Theorem 3.1}. Let $(X,Y)$ be a bipartition of bipartite graph $K_{m_{1},m_{2}}$(with $2 \leq m_{1} \leq m_{2}$), where $X=\{x_1,x_2,\cdots, x_{m_{1}}\},$ $Y=\{y_1,y_2,\cdots, y_{m_{2}}\}.$ If $S\subseteq X$ and $|S|\geq 2,$ then $\kappa_{K_{m_{1},m_{2}}}^*(S)=m_{2}.$ If $S\subseteq Y$ and $|S|\geq 2,$ then $\kappa_{K_{m_{1},m_{2}}}^*(S)=m_{1}.$

{\bf Proof.} Assume $S=\{x_1,x_2,\cdots,x_s\}\subseteq X.$ Let $T_j$ be an $S$-Steiner tree with $V(T_j)=S\cup \{y_{j}\},$ $E(T_j)=\{x_1y_j,x_2y_j,\cdots,x_sy_j\},$ where $1\leq j\leq m_{2}.$ By Theorem 2.1, $T_i(1\leq i\leq k)$ are $m_{2}$ CISSTs. So $\kappa_{K_{m_{1},m_{2}}}^*(S)\geq m_{2}.$ Since $S\subseteq X$ and $s\geq 2,$ by the definition of $S$-Steiner trees, every $S$-Steiner tree contains at least one vertex $y\in Y$ as an internal vertex. And by Theorem 2.1, the number of CISSTs in $K_{m_{1},m_{2}}$ is not more than $m_{2},$ that is, $\kappa_{K_{m_{1},m_{2}}}^*(S)\leq m_{2}.$ So $\kappa_{K_{m_{1},m_{2}}}^*(S)=m_{2}.$

Similarly, if $S\subseteq Y$ and $|S|\geq 2,$ then $\kappa_{K_{m_{1},m_{2}}}^*(S)=m_{1}.$
$\Box$

\begin{figure}[htbp]
  \centering
  \includegraphics[width=0.4\textwidth]{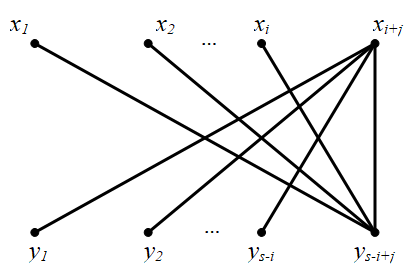}
  \caption{the $S_i$-Steiner tree $T_j$}
  \label{fig:figure4}
\end{figure}

Let $S_i=\{x_1,x_2,\cdots,x_i,y_1,y_2,\cdots,y_{s-i}\}.$ Assume $\min\{m_{1}-i,m_{2}-(s-i)\}\geq 1.$ Let $j$ be an integer with $1\leq j\leq \min\{m_{1}-i,m_{2}-(s-i)\},$ and we define an $S_i$-Steiner tree $T_j$ with $$V(T_j)=S_i\cup \{x_{i+j},y_{s-i+j}\},$$ $$E(T_j)=\{y_{s-i+j}x_1,y_{s-i+j}x_2,\cdots,y_{s-i+j}x_i\}\cup\{x_{i+j}y_1,x_{i+j}y_2, \cdots,x_{i+j}y_{s-i}\}\cup\{x_{i+j}y_{s-i+j}\},$$ where $x_{i+j}\in X\setminus S_i, y_{s-i+j}\in Y\setminus S_i$, see Figure 4. For $1\leq j\leq \min\{m_{1}-i,m_{2}-(s-i)\},$ the $S_i$-Steiner tree $T_j$ is said to be $I$-type, and the set of all $I$-type trees is denoted as ${\mathcal{A}}_{1}.$

Assume $1\leq i\leq \min\{m_{1}-1,s-1\}$ and $m_{2}-(s-i)<m_{1}-i.$ Then $|{\mathcal{A}}_{1}|=m_{2}-(s-i),$ and denote $a_1=|{\mathcal{A}}_{1}|.$ Let $k$ be an integer with $1\leq k\leq \min\{s-i,m_{1}-i-a_1\},$ and we define an $S_i$-Steiner tree $T_k'$ with $$V(T_k')=S_i\cup \{x_{i+a_1+k}\},$$ $$E(T_k')=\{x_{i+a_1+k}y_1,x_{i+a_1+k}y_2,\cdots,x_{i+a_1+k}y_{s-i}\}\cup\{y_k x_1,y_k x_2,\cdots,y_k x_i\},$$ see Figure 5. For $1\leq k\leq \min\{s-i,m_{1}-i-a_1\},$ the $S_i$-Steiner tree $T_k$ is said to be $I_X$-type, and the set of all $I_X$-type trees is denoted as ${\mathcal{A}}_{2}^{1}.$ On the other hand, consider a complete bipartite subgraph $K_{m_{1}-a_1,s-i}(X', Y')$ with $X'=\{x_1,x_2,\cdots,x_i,x_{i+a_1+1},x_{i+a_1+2},\cdots,$ $x_{m_{1}}\},$ $Y'=\{y_1,y_2,\cdots,y_{s-i}\}.$ The set of all completely independent spanning trees in $K_{m_{1}-a_1,s-i}(X', Y')$ is denoted as ${\mathcal{A'}}_{2}^{2}$. Since $S_{i}\subseteq V(K_{m_{1}-a_1,s-i}(X', Y'))$, by Theorem 2.2, there exists a completely independent $S_i$-Steiner tree set ${\mathcal{A}}_{2}^{2}$ such that $|{\mathcal{A}}_{2}^{2}|=|{\mathcal{A'}}_{2}^{2}|$.

\begin{figure}[htbp]
  \centering
  \includegraphics[width=0.4\textwidth]{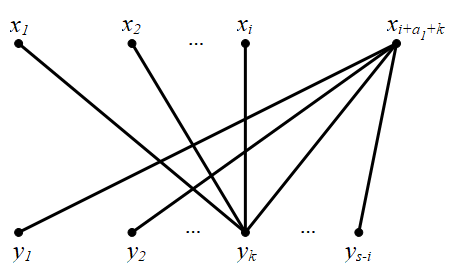}
  \caption{the $S_i$-Steiner tree $T_k'$}
  \label{fig:figure5}
\end{figure}

\begin{figure}[htbp]
  \centering
  \includegraphics[width=0.4\textwidth]{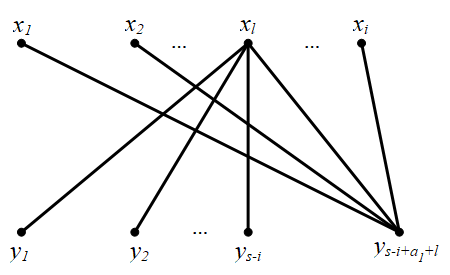}
  \caption{the $S_i$-Steiner tree $T_l''$}
  \label{fig:figure6}
\end{figure}

Assume $i\geq 1$ and $m_{2}-(s-i)>m_{1}-i.$ Then $a_1=|{\mathcal{A}}_{1}|=m_{1}-i.$ Let $l$ be an integer with $1\leq l\leq \min\{i,m_{2}-(s-i)-a_1\}$, and we define an $S_i$-Steiner tree $T_l''$ with $$V(T_l'')=S_i\cup \{y_{s-i+a_1+l}\},$$ $$E(T_l'')=\{y_{s-i+a_1+l}x_1,y_{s-i+a_1+l}x_2, \cdots,y_{s-i+a_1+l}x_i\}\cup\{x_l y_1, x_l y_2,\cdots, x_l y_{s-i}\},$$ see Figure 6. For $1\leq l\leq \min\{i,m_{2}-(s-i)-a_1\},$ the $S_i$-Steiner tree $T_l''$ is said to be $I_Y$-type, and the set of all $I_Y$-type trees is denoted as ${\mathcal{A}}_{3}^{1}.$ On the other hand, consider a complete bipartite subgraph $K_{i,m_{2}-a_1}(X'', Y'')$ with $X''=\{x_1,x_2,\cdots,x_i\}, Y''=\{y_1,y_2,\cdots,y_{s-i}, y_{s-i+a_1+1},y_{s-i+a_1+2},\cdots,y_{m_{2}}\}.$ The set of all completely independent spanning trees in $K_{i,m_{2}-a_1}(X'', Y'')$ is denoted as ${\mathcal{A'}}_{3}^{2}.$ Since $S_{i}\subseteq V(K_{i,m_{2}-a_1}(X'',Y''))$, by Theorem 2.2, there exists a completely independent $S_i$-Steiner tree set ${\mathcal{A}}_{3}^{2}$ such that $|{\mathcal{A}}_{3}^{2}|=|{\mathcal{A'}}_{3}^{2}|$.

Based on the above construction, we can easily obtain the following lemma.

{\bf Lemma 3.2}. The $S_i$-Steiner trees in ${\mathcal{A}}_{1}$ are completely independent, and $a_1=\min\{m_{1}-i,m_{2}-(s-i)\}.$ Furthermore, if $a_1=m_{2}-(s-i),$ then $|{\mathcal{A}}_{2}^{1}|= \min\{s-i,m_{1}-i-a_1\},$ and the $S_i$-Steiner trees in ${\mathcal{A}}_{1}\cup {\mathcal{A}}_{2}^{1}$ are completely independent; if $a_1=m_{1}-i,$then $|{\mathcal{A}}_{3}^{1}|= \min\{i,m_{2}-(s-i)-a_1\},$ and the $S_i$-Steiner trees in ${\mathcal{A}}_{1}\cup {\mathcal{A}}_{3}^{1}$ are completely independent.

{\bf Lemma 3.3}\cite{Pai}. Let $K_{m_{1},m_{2}}$ be a complete bipartite graph, and $m_{1}\leq m_{2}.$ Then $K_{m_{1},m_{2}}$ contains $\lfloor\frac{m_{1}}{2}\rfloor$ CISTs.

Lemma 3.3 leads to the following conclusion.

{\bf Corollary 3.4}. If $a_1=m_{2}-(s-i),$ then $|{\mathcal{A}}_{2}^{2}|=\max\{\min\{\lfloor\frac{m_{1}-a_1}{2}\rfloor,\lfloor\frac{s-i}{2}\rfloor\},1\},$ and the $S_i$-Steiner trees in ${\mathcal{A}}_{1}\cup {\mathcal{A}}_{2}^{2}$ are completely independent; if $a_1=m_{1}-i,$ then $|{\mathcal{A}}_{3}^{2}|=\max\{\min\{\lfloor\frac{i}{2}\rfloor,\lfloor\frac{m_{2}-a_1}{2}\rfloor\},1\},$ and the $S_i$-Steiner trees in ${\mathcal{A}}_{1}\cup {\mathcal{A}}_{3}^{2}$ are completely independent.

{\bf Theorem 3.5}. Let $(X,Y)$ be a bipartition of the bipartite graph $K_{m_{1},m_{2}}$(with $2 \leq m_{1} \leq m_{2}$), where $X=\{x_1,x_2,\cdots, x_{m_{1}}\},Y=\{y_1,y_2,\cdots, y_{m_{2}}\}.$ And let $S_i=\{x_1,x_2,\cdots, x_i,y_1,y_2,\cdots, y_{s-i}\}$ be a subset of $s$ vertices in $K_{m_{1},m_{2}},$ where $\max\{1, s-m_{2}\}\leq i \leq \min\{m_{1},s-1\}.$ If $s\leq m_{2}-m_{1}+2,$ then $\kappa_{K_{m_{1},m_{2}}}^*(S_i)= m_{1}.$

{\bf Proof.} Since $ i \leq s-1,$ it follows that $S_{i}\bigcap Y\neq \emptyset.$ Note that for every vertex $w \in S_{i}\bigcap Y$, $d_{K_{m_{1},m_{2}}}(w) = m_{1}.$ And by any two completely independent $S_i$-Steiner trees are edge-disjoint, we can conclude that the number of completely independent $S_i$-Steiner trees in $K_{m_{1},m_{2}}$ is not more than $m_{1},$ that is, $\kappa_{K_{m_{1},m_{2}}}^*(S_i)\leq m_{1}.$

Next, we shall show $\kappa_{K_{m_{1},m_{2}}}^*(S_i)\geq m_{1}.$ Consider $i=1.$ Since $m_{1}\geq 2,$ we have $m_{1}-i\geq 1.$ By $s\leq m_{2}-m_{1}+2,$ it follows that $m_{2}-(s-i)=m_{2}-s+i\geq (m_{1}-2)+i=m_{1}-i.$ So $\min\{m_{1}-i,m_{2}-(s-i)\}=m_{1}-i\geq 1.$ By Lemma 3.2, $|{\mathcal{A}}_{1}|=m_{1}-1.$ Clearly, $K_{m_{1},m_{2}}[S_1]$ is an $S_i$-Steiner tree, denoted by $T_0$. By Theorem 2.1, $\{T_0\}\cup {\mathcal{A}}_{1}$ contains $m_{1}$ completely independent $S_i$-Steiner trees in $K_{m_{1},m_{2}}.$ So $\kappa_{K_{m_{1},m_{2}}}^*(S_i)\geq m_{1}.$

Assume $2\leq i\leq m_{1}-1.$ It follows that $m_{1}-i\geq 1.$ By $s\leq m_{2}-m_{1}+2,$ It follows that $m_{2}-(s-i)=m_{2}-s+i\geq (m_{1}-2)+i\geq m_{1}-i.$ So $\min\{m_{1}-i,m_{2}-(s-i)\}=m_{1}-i\geq 1.$ By Lemma 3.2, $|{\mathcal{A}}_{1}|=m_{1}-i.$ In addition, by $s\leq m_{2}-m_{1}+2,$ we also deduce that $m_{2}-(s-i)-(m_{1}-i)=m_{2}-s-m_{1}+2i\geq i.$ By Lemma 3.2, $|{\mathcal{A}}_{3}^{1}|=\min\{i,m_{2}-(s-i)-(m_{1}-i)\}=i.$ Furthermore, ${\mathcal{A}}_{3}^{1}\cup {\mathcal{A}}_{1}$ contains $m_{1}$ completely independent $S_i$-Steiner trees in $K_{m_{1},m_{2}}.$ So $\kappa_{K_{m_{1},m_{2}}}^*(S_i)\geq m_{1}.$

Consider $i=m_{1}.$ By Lemma 3.2, $|{\mathcal{A}}_{1}|=0.$ By $s\leq m_{2}-m_{1}+2,$ we also deduce that $m_{2}-(s-i)=m_{2}-s+i\geq (m_{1}-2)+i\geq i.$ So $\min\{i,m_{2}-(s-i)\}\geq i\geq 2.$ By Lemma 3.2, $|{\mathcal{A}}_{3}^{1}|=\min\{i,m_{2}-(s-i)\}=i=m_{1}.$ By the definition of completely independent $S_i$-Steiner trees, ${\mathcal{A}}_{3}^{1}$ contains $m_{1}$ completely independent $S_i$-Steiner trees in $K_{m_{1},m_{2}}.$ So $\kappa_{K_{m_{1},m_{2}}}^*(S_i)\geq m_{1}.$ The proof is complete.
$\Box$

{\bf Theorem 3.6}. Let $(X,Y)$ be a bipartition of the bipartite graph $K_{m_{1},m_{2}}$(with $2 \leq m_{1} \leq m_{2}$), where $X=\{x_1,x_2,\cdots, x_{m_{1}}\},Y=\{y_1,y_2,\cdots, y_{m_{2}}\}.$ And let $S_i=\{x_1,x_2,\cdots, x_i,y_1,y_2,\cdots, y_{s-i}\}$ be a subset of $s$ vertices in $K_{m_{1},m_{2}},$ where $\max\{1, s-m_{2}\}\leq i \leq \min\{m_{1},s-1\}$ and $s\geq m_{2}-m_{1}+3.$ If $2\leq 2i\leq m_{1}+s-m_{2},$ then
\begin{equation*}
\kappa_{K_{m_{1},m_{2}}}^*(S_i)\geq
\begin{cases}
m_{1}-i,& \text{ $ 2\leq 2i\leq \frac{2(m_{1}+s-m_{2})}{3}; $}\\
m_{2}-(s-i)+\lfloor\frac{m_{1}+s-m_{2}-i}{2}\rfloor, & \text{ $ \frac{2(m_{1}+s-m_{2})}{3}< 2i\leq m_{1}+s-m_{2}. $}
\end{cases}
\end{equation*}

{\bf Proof.} Since $2i\leq m_{1}+s-m_{2},$ we have $m_{2}-(s-i)\leq m_{1}-i.$ So $\min\{m_{1}-i,m_{2}-(s-i)\}=m_{2}-(s-i).$ By Lemma 3.2, $a_1=|{\mathcal{A}}_{1}|=m_{2}-(s-i).$ It follows that $ m_{1}-a_1=m_{1}+s-m_{2}-i=s-i+(m_{1}-m_{2})\leq s-i.$ And by $2\leq 2i\leq m_{1}+s-m_{2},$ we have $m_{1}-a_1=m_{1}+s-m_{2}-i\geq i\geq 1.$ If $m_{1}-a_1=1,$ then $i=1$ and $m_{1}+s-m_{2}=2,$ contradicting $s\geq m_{2}-m_{1}+3.$ So assume $m_{1}-a_1\geq 2.$ By Corollary 3.4, we have $|{\mathcal{A}}_{2}^{2}|=\max\{\min\{\lfloor\frac{m_{1}-a_1}{2}\rfloor,\lfloor\frac{s-i}{2}\rfloor\},1\}=\lfloor\frac{m_{1}-a_1}{2}\rfloor=\lfloor\frac{m_{1}+s-m_{2}-i}{2}\rfloor.$ By $a_1=m_{2}-(s-i),$ we have $m_{1}-i-a_1=s-i+(m_{1}-m_{2}-i)<s-i.$ By Lemma 3.2, we have $|{\mathcal{A}}_{2}^{1}|= \min\{s-i,m_{1}-i-a_1\}=m_{1}-i-a_1=m_{1}+s-m_{2}-2i.$

By Lemma 3.2 and Corollary 3.4, $S_i$-Steiner trees in ${\mathcal{A}}_{1}\cup {\mathcal{A}}_{2}^{1}$ are completely independent, and the $S_i$-Steiner trees in ${\mathcal{A}}_{1}\cup {\mathcal{A}}_{2}^{2}$ are also completely independent. Thus $\kappa_{K_{m_{1},m_{2}}}^*(S_i)\geq \max\{|{\mathcal{A}}_{1}|+|{\mathcal{A}}_{2}^{1}|, |{\mathcal{A}}_{1}|+|{\mathcal{A}}_{2}^{2}|\}.$

If $2i\leq \frac{2(m_{1}+s-m_{2})}{3},$ then $m_{1}+s-m_{2}-2i\geq \lfloor\frac{m_{1}+s-m_{2}-i}{2}\rfloor,$ hence $\kappa_{K_{m_{1},m_{2}}}^*(S_i)\geq |{\mathcal{A}}_{1}|+|{\mathcal{A}}_{2}^{1}|=m_{1}-i.$ If $\frac{2(m_{1}+s-m_{2})}{3}<2i\leq m_{1}+s-m_{2},$ then $m_{1}+s-m_{2}-2i\leq \lfloor\frac{m_{1}+s-m_{2}-i}{2}\rfloor,$ $\kappa_{K_{m_{1},m_{2}}}^*(S_i)\geq |{\mathcal{A}}_{1}|+|{\mathcal{A}}_{2}^{2}|=m_{2}-(s-i)+\lfloor\frac{m_{1}+s-m_{2}-i}{2}\rfloor.$ The proof is complete.
$\Box$

{\bf Theorem 3.7}. Let $(X,Y)$ be a bipartition of the bipartite graph $K_{m_{1},m_{2}}$(with $2 \leq m_{1} \leq m_{2}$), where $X=\{x_1,x_2,\cdots, x_{m_{1}}\},Y=\{y_1,y_2,\cdots, y_{m_{2}}\}.$ And let $S_i=\{x_1,x_2,\cdots, x_i,y_1,y_2,\cdots, y_{s-i}\}$ be a subset of $s$ vertices in $K_{m_{1},m_{2}},$ where $\max\{1, s-m_{2}\}\leq i \leq \min\{m_{1},s-1\}$ and $s\geq m_{2}-m_{1}+3.$ If $2i>m_{1}+s-m_{2},$ then
\begin{equation*}
\kappa_{K_{m_{1},m_{2}}}^*(S_i)\geq
\begin{cases}
m_{1}-i+\lfloor\frac{i}{2}\rfloor,& \text{ $ m_{1}+s-m_{2}< 2i\leq \frac{4(m_{1}+s-m_{2})}{3}; $}\\
m_{2}-s+i, & \text{ $\frac{4(m_{1}+s-m_{2})}{3}< 2i\leq 2(m_{1}+s-m_{2});$}\\
m_{1}, & \text{ $2i> 2(m_{1}+s-m_{2}).$}\\
\end{cases}
\end{equation*}

{\bf Proof.} Since $2i> m_{1}+s-m_{2},$ we have $m_{2}-(s-i)> m_{1}-i.$ It follows that $\min\{m_{1}-i,m_{2}-(s-i)\}=m_{1}-i.$ By Lemma 3.2, $a_1=|{\mathcal{A}}_{1}|=\min\{m_{1}-i,m_{2}-(s-i)\}=m_{1}-i.$ It follows that $m_{2}-a_1=m_{2}-m_{1}+i\geq i.$ Since $s\geq m_{2}-m_{1}+3$ and $2i>m_{1}+s-m_{2},$ we have $2i>3$ and hence $i\geq 2.$ By Corollary 3.4, $|{\mathcal{A}}_{3}^{2}|=\max\{\min\{\lfloor\frac{i}{2}\rfloor,\lfloor\frac{m_{2}-a_1}{2}\rfloor\},1\}=\lfloor\frac{i}{2}\rfloor.$ By Lemma 3.2 and Corollary 3.4, the $S_i$-Steiner trees in ${\mathcal{A}}_{1}\cup {\mathcal{A}}_{3}^{1}$ are completely independent, and the $S_i$-Steiner trees in ${\mathcal{A}}_{1}\cup {\mathcal{A}}_{3}^{2}$ are also completely independent. Thus $\kappa_{K_{m_{1},m_{2}}}^*(S_i)\geq \max\{|{\mathcal{A}}_{1}|+|{\mathcal{A}}_{3}^{1}|, |{\mathcal{A}}_{1}|+|{\mathcal{A}}_{3}^{2}|\}.$ We discuss the following two cases.

{\bf Case 1}. $m_{1}+s-m_{2}<2i\leq 2(m_{1}+s-m_{2}).$

Then $1\leq m_{2}-m_{1}-s+2i \leq i.$ By Lemma 3.2, we have $|{\mathcal{A}}_{3}^{1}|= \min\{i,m_{2}-(s-i)-a_1\}=\min\{i,m_{2}-m_{1}-s+2i\}=m_{2}-m_{1}-s+2i\geq 1.$

{\bf Case 1.1}. $m_{1}+s-m_{2}<2i\leq \frac{4(m_{1}+s-m_{2})}{3}.$

In this case, we have $m_{2}-s-m_{1}+2i\leq \lfloor\frac{i}{2}\rfloor,$ that is, $|{\mathcal{A}}_{3}^{1}|\leq |{\mathcal{A}}_{3}^{2}|.$ It follows that $\kappa_{K_{m_{1},m_{2}}}^*(S_i)\geq \max\{|{\mathcal{A}}_{1}|+|{\mathcal{A}}_{3}^{1}|, |{\mathcal{A}}_{1}|+|{\mathcal{A}}_{3}^{2}|\}=|{\mathcal{A}}_{1}|+|{\mathcal{A}}_{3}^{2}|=m_{1}-i+\lfloor\frac{i}{2}\rfloor.$

{\bf Case 1.2}. $\frac{4(m_{1}+s-m_{2})}{3}<2i\leq 2(m_{1}+s-m_{2}).$

In this case, we have $m_{2}-s-m_{1}+2i> \lfloor\frac{i}{2}\rfloor,$ that is, $|{\mathcal{A}}_{3}^{1}|\geq |{\mathcal{A}}_{3}^{2}|.$ It follows that $\kappa_{K_{m_{1},m_{2}}}^*(S_i)\geq \max\{|{\mathcal{A}}_{1}|+|{\mathcal{A}}_{3}^{1}|, |{\mathcal{A}}_{1}|+|{\mathcal{A}}_{3}^{2}|\}=|{\mathcal{A}}_{1}|+|{\mathcal{A}}_{3}^{1}|=m_{1}-i+m_{2}-m_{1}-s+2i=m_{2}-s+i.$

{\bf Case 2}. $2i> 2(m_{1}+s-m_{2}).$

By $a_1=m_{1}-i$, we have $m_{2}-(s-i)-a_1=m_{2}-s-m_{1}+2i>i.$ By Lemma 3.2, we have $|{\mathcal{A}}_{3}^{1}|= \min\{i,m_{2}-(s-i)-a_1\}=\min\{i,m_{2}-m_{1}-s+2i\}=i.$ Note that $|{\mathcal{A}}_{3}^{2}|=\lfloor\frac{i}{2}\rfloor.$ So $\kappa_{K_{m_{1},m_{2}}}^*(S_i)\geq \max\{|{\mathcal{A}}_{1}|+|{\mathcal{A}}_{3}^{1}|, |{\mathcal{A}}_{1}|+|{\mathcal{A}}_{3}^{2}|\}=|{\mathcal{A}}_{1}|+|{\mathcal{A}}_{3}^{1}|=m_{1}-i+i=m_{1}.$
The proof has been completed. $\Box$

{\bf Corollary 3.8}. Let $(X,Y)$ be a bipartition of the bipartite graph $K_{m_{1},m_{2}}$(with $2 \leq m_{1} \leq m_{2}$), where $|X|=m_{1}, |Y|=m_{2},$ and let $s$ be an integer with $s\geq 2.$ Then we have

(1) If $s\leq m_{2}-m_{1}+2,$ then $\kappa_s^*(K_{m_{1},m_{2}})=m_{1}$ except $S\subseteq X.$

(2) If $S\subseteq X,$ then $\kappa_s^*(K_{m_{1},m_{2}})=m_{2}$ for any $2\leq s\leq m_{1}.$

(3) If $s\geq m_{2}-m_{1}+3,$ then $\kappa_s^*(K_{m_{1},m_{2}})\geq m_{1}-\frac{m_{1}+s-m_{2}+2}{3}.$

{\bf Proof.} By Theorems 3.1 and 3.5, Statements (1) and (2) are true.

For convenience, let $z=m_{1}+s-m_{2}.$ By Theorems 3.6 and 3.7, we can construct the function
\begin{equation*}
f(i)=
\begin{cases}
m_{1}-i,& \text{ $ 2\leq 2i\leq \frac{2z}{3}; $}\\
m_{2}-(s-i)+\lfloor\frac{m_{1}+s-m_{2}-i}{2}\rfloor, & \text{ $ \frac{2z}{3}< 2i\leq z; $}\\
m_{1}-i+\lfloor\frac{i}{2}\rfloor,& \text{ $ z< 2i\leq \frac{4z}{3}; $}\\
m_{2}-s+i, & \text{ $\frac{4z}{3}< 2i\leq 2z;$}\\
m_{1}, & \text{ $2i> 2z.$}\\
\end{cases}
\end{equation*}
According to the monotonicity property of this function, we have
\begin{equation*}
\min f(i) = \min \left\{
    f(\left\lfloor \frac{z}{3} \right\rfloor),
    f(\left\lceil \frac{z}{3} \right\rceil),
     f(\left\lfloor \frac{2z}{3} \right\rfloor),
     f(\left\lceil \frac{2z}{3} \right\rceil)
\right\}.
\end{equation*}
By the definition of this function, we derive the following inequalities

$f(\left\lfloor \frac{z}{3} \right\rfloor) \geq m_{1}- \frac{m_{1}+s-m_{2}}{3} > m_{1} - \frac{m_{1}+s-m_{2}+2}{3};$

$f(\left\lceil \frac{z}{3} \right\rceil)\geq m_{2}-\left( s-\frac{m_{1}+s-m_{2}}{3} \right)+\lfloor\frac{m_{1}+s-m_{2}-\frac{m_{1}+s-m_{2}}{3}}{2}\rfloor=m_{1}-\frac{2(m_{1}+s-m_{2})}{3}+\lfloor\frac{m_{1}+s-m_{2}}{3}\rfloor\geq m_{1} - \frac{m_{1}+s-m_{2}+2}{3};$

$f(\left\lfloor \frac{2z}{3} \right\rfloor) \geq m_{1} - \left( \frac{2(m_{1}+s-m_{2})}{3} \right)+\lfloor\frac{\frac{2(m_{1}+s-m_{2})}{3}}{2}\rfloor=m_{1}-\frac{2(m_{1}+s-m_{2})}{3}+\lfloor\frac{m_{1}+s-m_{2}}{3}\rfloor\geq m_{1} - \frac{m_{1}+s-m_{2}+2}{3};$

$f(\left\lceil \frac{2z}{3} \right\rceil)\geq m_{2}- s + \frac{2(m_{1}+s-m_{2})}{3} =m_{1} - \frac{m_{1}+s-m_{2}}{3} > m_{1} - \frac{m_{1}+s-m_{2}+2}{3}.$

So, $\min f(i)\geq m_{1}-\frac{m_{1}+s-m_{2}+2}{3}$, it suggests $\kappa_s^*(K_{m_{1},m_{2}})\geq m_{1}-\frac{m_{1}+s-m_{2}+2}{3}.$$\Box$

{\bf Lemma 3.9}. Let $(X,Y)$ be a bipartition of the bipartite graph $K_{5,6}$, where $X=\{x_1,x_2,\cdots, x_5\},$ $Y=\{y_1,y_2,\cdots, y_6\}.$ And let $S=\{x_1,x_2,y_1,y_2\}.$ Then $\kappa_{K_{5,6}}^*(S)=4.$

{\bf Proof.} Since $|S\cap X|=2$ and $|S\cap Y|=2,$ it follows that in any $S$-Steiner tree $T,$ there must exist an internal vertex $x\in X$ that connects the vertices in $Y \cap S,$ and there must also exist an internal vertex $y\in Y$ that connects the vertices in $X\cap S.$ By Theorem 2.1, every vertex in $X$ is an internal vertex of at most one of the CISSTs. Thus, $\kappa_{K_{5,6}}^*(S)\leq |X|=5.$

Suppose $\kappa_{K_{5,6}}^*(S)=5,$ and say $T_1,T_2,\cdots,T_5$ are $5$ CISSTs. Note that every vertex in $K_{5,6}$ is at most an internal vertex of $T_{i}(1\leq i\leq5).$ And by the symmetry of $K_{5,6},$ without loss of generality, assume $x_i,y_i$ are two internal vertices in $T_i$ for $i=1,2,\cdots,5.$ Since $x_i$ is the sole internal vertex of $V(T_i)\cap X,$ we have $x_i y_1,x_i y_2\in E(T_i)$ for $i=1,2,\cdots,5.$ It follows that $x_1 y_1\in E(T_1).$ Consider the $S$-Steiner tree $T_1.$ Since $y_1$ is the internal vertex of $T_1,$ there is a vertex $x_t,t\neq 1$ such that $x_t y_1\in E(T_1).$ On the other hand, since $x_i y_1\in E(T_i)$ for $i=1,2,\cdots,5,$ we have $x_t y_1\in E(T_t).$ It follows that $x_t y_1\in E(T_1)\cap E(T_t),$ a contradiction. So $\kappa_{K_{5,6}}^*(S)<5.$

Moreover, let $m_{1}=5,m_{2}=6,s=4,i=2.$ Then $2i=\frac{4(m_{1}+s-m_{2})}{3}.$ By Theorem 3.7, we can deduce that $\kappa_{K_{5,6}}^*(S)\geq 4.$ Therefore, $\kappa_{K_{5,6}}^*(S)=4.$ The proof is complete. $\Box$

{\bf Remark 3.10}. Let $m_{1}=5,m_{2}=6,s=4,i=2.$ By Lemma 3.9, we have $\kappa_{K_{m_{1},m_{2}}}^*(S)=m_{1}-\lfloor\frac{m_{1}+s-m_{2}+2}{3}\rfloor.$ In this sense, the lower bound of Corollary 3.8(3) is tighter.

\section{Conclusions}

In this paper, we introduce the novel concepts of completely independent $S$-Steiner trees (CISSTs in short) and generalized $k^*$-connectivity in graph theory. We establish a theoretical characterization theorem for CISSTs, and investigate the properties of CISSTs and generalized $k^*$-connectivity. Building upon this theoretical framework, we investigate the CISSTs in complete graphs and complete bipartite graphs. Furthermore, we determine the generalized $k^*$-connectivity for complete graphs and give a tight lower bound of generalized $k^*$-connectivity for complete bipartite graphs.

\section*{Acknowledgments}
This work was financially supported by Basic Research Projects of Shanxi Province (202203021221154), and the Natural Science Foundation of Shanxi Province (202203021221037).

\section*{Declaration of competing interest}
The authors declare that they have no known competing financial interests or personal relationships that could have appeared to influence the work reported in this paper.


\begin{thebibliography}{99}

\bibitem{DuH} Du D, Hu X. Steiner Tree Problems in Computer Communication Networks. World Scientific; 2008.

\bibitem{Whitney2} Whitney H. Congruent graphs and the connectivity of graphs. Amer J Math. 1932;54:150-168.
\bibitem{Hager1985} Hager M. Pendant tree-connectivity. J Combin Theory Ser B. 1985;38(2):179-189.
\bibitem{Lin} Lin SW, Zhang QH. The generalized 4-connectivity of hypercubes. Discrete Appl Math. 2017;220:60-67.
\bibitem{Zhao'} Zhao SL, Hao RX. The generalized 4-connectivity of exchanged hypercubes. Appl Math Comput. 2019;347:342-353.
\bibitem{Zhao2} Zhao SL, Hao RX, Cheng E. Two kinds of generalized connectivity of dual cubes. Discrete Appl Math. 2019;257:306-316.
\bibitem{Li2020} Li CF, Lin SW, Li SJ. The 4-set tree connectivity of $(n,k)$-star networks. Theor Comput Sci. 2020;844:81-86.
\bibitem{Zhao21} Zhao SL, Hao RX, Wu J. The generalized 4-connectivity of hierarchical cubic networks. Discrete Appl Math. 2021;289:194-206.
\bibitem{Chartrand} Chartrand G, Okamoto F, Zhang P. Rainbow trees in graphs and generalized connectivity. Networks. 2010;55(4):360-367.
\bibitem{Li2010'} Li SS, Li W, Li XL. The generalized connectivity of complete bipartite graphs. Ars Combin. 2010;104:65-79.
\bibitem{Li2014} Li SS, Li W, Li XL. The generalized connectivity of complete equipartition 3-partite graphs. Bull Malays Math Sci Soc. 2014;37(1):103-121.
\bibitem{Li2012} Li SS, Li XL. Note on the hardness of generalized connectivity. J Combin Optimization. 2012;24(3):389-396.
\bibitem{Li2010} Li SS, Li XL, Zhou WL. Sharp bounds for the generalized connectivity $\kappa_3(G)$. Discrete Math. 2010;310(13-14):2147-2163.
\bibitem{Hasunuma2} Hasunuma T. Completely independent spanning trees in the underlying graph of a line digraph. Discrete Math. 2001;234(1-3):149-157.
\bibitem{Chen L} Chen L, Li XL, Liu MM, et al. Further hardness results on the generalized (edge-)connectivity of graphs. arXiv:1304.6153 [math.CO]. 2013.
\bibitem{Pai} Pai KJ, Tang SM, Chang JM, et al. Completely independent spanning trees on complete graphs, complete bipartite graphs and complete tripartite graphs. Adv Intell Syst Appl. 2013;1(1):107-113.

\end{thebibliography}
\end{document}